\documentclass[12pt, reqno]{amsart} \usepackage{amsmath,
amssymb}

\numberwithin{equation}{section} 

\newcommand\A{{\mathbf A}} \newcommand\bP{{\mathbf P}}
 \newcommand\F{{\mathbf F}}
\newcommand\G{{\mathbf G}} \newcommand\Z{{\mathbf Z}}
 
 \newcommand\Q{{\mathbf Q}}
\newcommand\N{{\mathbf N}}

\newcommand\calS{{\mathsf S}} 
\newcommand\calO{{\mathcal O}}
\newcommand{\clO}{\mathcal{O}} 
\newcommand{\lb}{[} \newcommand{\rb}{]} \newcommand{\brk}[1]{\lb#1\rb}
\newcommand{\spn}[1]{\langle #1\rangle} \newcommand{\num}[1]{\# #1}
\newcommand{\cmapsto}{\rightsquigarrow}
\newcommand{\primaries}{\mathcal{Q}} 
\newcommand{\ring}{{\mathsf R}}

\newcommand{\dom}{\operatorname {dom}}
\newcommand{\coker}{\operatorname {coker}}
\newcommand\Gr{\operatorname {Gr}} \newcommand\Gl{\operatorname {GL}}
\newcommand\GL{\operatorname {GL}}  
\newcommand\Sym{\operatorname {Sym}}
\newcommand\Spec{\operatorname{Spec}}  
\newcommand\Hom{\operatorname {Hom}}  \newcommand\Fun{\operatorname {Fun}}

\newcommand\union{\cup}

\newcommand\beea{\begin{eqnarray*}} \newcommand\eeea{\end{eqnarray*}}
\newcommand\eea{\end{eqnarray*}} \newcommand\beean{\begin{eqnarray}}
\newcommand\eeean{\end{eqnarray}} \newcommand\beqn{\begin{equation}}
\newcommand\beq{\begin{equation}} 
\newcommand\ee{\end{equation}} \newcommand\eeqn{\end{equation}}

\newcommand\Mot{\operatorname{\sf Mot}} 
\newcommand\Graphs{\operatorname{\sf Graphs}}
\newcommand\Matroids{\operatorname{\sf Matroids}}

\newtheorem{theorem}{Theorem}[section]
\newtheorem{lemma}[theorem]{Lemma}
\newtheorem{proposition}[theorem]{Proposition}
\newtheorem{corollary}[theorem]{Corollary}
\newtheorem{definition-proposition}[theorem]{Definition-Proposition}

\theoremstyle{definition} \newtheorem{definition}[theorem]{Definition}
\newtheorem{example}[theorem]{Example}
\newtheorem{conjecture}[theorem]{Conjecture}

\theoremstyle{remark} \newtheorem{remark}[theorem]{Remark}

\begin{document}
\title[Matroids, motives and a conjecture of Kontsevich]{Matroids,
motives and a conjecture of Kontsevich} \author{Prakash Belkale}
\address{Department of Mathematics\\ University of Utah \\ 155 South
1400 East, JWB 233 \\ Salt Lake City, Utah 84112-0090}
\email{belkale@math.utah.edu} \author{Patrick Brosnan}
\address{Department of Mathematics\\ University of California,
Irvine\\ Irvine, California 92697}
% Max-Planck-Institute f\"{u}r Mathematik\\ Vivatsgasse 7\\
% D-53111 Bonn\\ Germany}
\email{pbrosnan@math.uci.edu}
\begin{abstract} 
Let $G$ be a finite connected graph with $\# E$ edges and with first
betti number $b_1$.  The Kirchhoff polynomial $P_G$ is a certain
homogeneous polynomial of degree $b_1$ in $\# E$ variables.  These
polynomials appear in the study of electrical circuits and in the
evaluation of Feynman amplitudes.  Motivated by work of D.~Kreimer and
D.~J.~Broadhurst associating multiple zeta values to certain Feynman
integrals, Kontsevich conjectured that the number of zeros of $P_G$
over the field with $q$ elements is a polynomial function of $q$.  We
show that this conjecture is false by relating the schemes $V(P_G)$ to
the representation spaces of matroids.  Moreover, using Mn\"ev's
universality theorem, we show that the schemes $V(P_G)$ essentially
generate all arithmetic of schemes over $\Z$.
\end{abstract}

\maketitle
\bigskip

\section[intro]{Introduction}

\subsection{Kontsevich's Conjecture}  Let $G$ be a finite graph with
vertex set $V= V(G)$ and edge set $E= E(G)$.  Considering such graphs
as finite CW-complexes, the betti numbers $b_0(G)$ and $b_1(G)$ are
both defined.  Recall that a graph $T$ is called a tree if $b_0(T)=1$
and $b_1(T)=0$.  A subgraph $T\subset G$ is called a spanning tree if
$T$ is a tree and $V(T)=V(G)$.  In a connected graph, a tree is a
spanning tree if and only if it is maximal.
 
For each edge $e$, let $x_e$ denote a formal variable.  Consider the
polynomial
\begin{equation}\label{eq:KontsevichPolynomial}
P_{G} = \sum_{T} \prod_{e\not\in T} x_e
\end{equation}
where the sum runs through all spanning trees of $G$.  If $G$ is not
connected, $P_G=0$ because the sum is empty.  Otherwise, $P_G$ is a
homogeneous polynomial of degree $b_1(G)$.

Let $V(P_G)$ denote the scheme of zeros of $P_G$ over $\Z$ --- a
hypersurface in $\A^E$.  Let $Y_G$ denote the complement of $V(P_G)$
in $\A^E$.  Motivated by computer calculations of the counterterms
appearing in the renormalization of Feynman integrals
\cite{kreimerbroadhurst}, Kontsevich speculated that the periods of
$Y_G$ are multiple zeta values (MZVs).  (See~\cite{kreimer} for a
discussion of MZVs and their relationship to certain Feynman
amplitudes).  Under this assumption on the periods, it is natural
to expect that the zeta functions associated to the $Y_G$ are the zeta
functions of motives of mixed Tate type \cite{ZagierZeta}.
  
Based on this idea, M.~Kontsevich made a conjecture about the number
of points of $Y_G$ over a finite field \cite{KGelfand}.
To describe
this conjecture, we first make a notational convention: For any scheme
$X$ of finite type over $\Z$, let $[X]$ denote the function $q\mapsto
\# X(\F_q)$.  Thus $\brk{X}$ is a function from the set $\primaries$
of prime powers to $\Z$.  We say that $X$ is {\em polynomially
countable} if $[X]$ is a polynomial in $\Z[q]$.

\begin{conjecture}[Kontsevich]\label{conjK}
For all graphs $G$, $[Y_G]\in\Z[q]$.
\end{conjecture}

Since $[V(P_G)] + [Y_G] = q^{\num{E}}$, this conjecture is equivalent
to the conjecture that $V(P_G)$ is polynomially countable.

Stembridge~\cite{stembridge} verified this conjecture for all graphs
with fewer than 12 edges.  For certain graphs it is relatively easy to
see that the conjecture holds.  For example, for $G$ a cycle of length
$n$, $V(P_G)$ is simply $\A^{n-1}$ and, thus, $[Y_G]= q^n - q^{n-1}$.

We will show, however, that Conjecture~\ref{conjK} is false.  In fact,
contrary to the extremely strong restrictions on the arithmetic nature
of the schemes $Y_G$ claimed by the conjecture, they are,
from an arithmetic point of view, the most general schemes possible.

\subsection{Motives and the Main Theorem}
To make this last statement precise we introduce some notation.  Let
$\Mot^+$ denote the group generated by all functions of the form $[X]$
for $X$ a scheme of finite type over $\Z$.  We think of $\Mot^+$ as a
coarse version of the ring of motives over $\Z$ 
\footnote{
In fact, all statements in this paper involving $\Mot^+$ do remain
valid in the finest possible setting.  That is we can
replace $\Mot^+$ with
the Grothendieck group defined by imposing the relation 
$[X]=[U]+[X-U]$ for $X$ a scheme of finite type over $\Z$, 
$U$ an open subset, and $[X]$ a formal symbol associated to $X$.}.
As $[X\times Y] =
[X][Y]$, $\Mot^+$ is a ring.  And, as $[\A^1] = q$, $\Mot^+$ is a
$\Z[q]$ module.  Let $\calS$ be the saturated multiplicative system in
$\Z[q]$ generated by the functions $q^n-q$ for $n>1$.  Set
$\Mot=\calS^{-1}\Mot^+$.  We remark that, since the functions in
$\calS$ are nonvanishing on $\primaries$,  
elements of $\Mot$ give everywhere-defined   
functions from $\primaries$ to $\Q$.

Let $\ring=\calS^{-1}\Z[q]$.  It is interesting to note that the
$K$-theory of $\ring$ turns up in questions about dynamical systems,
and, in a study of this $K$-theory \cite{grayson}, Grayson
showed that $\ring$ is a principal ideal domain.

Let $\Graphs$ denote the $\ring$-module generated by all functions of
the form $[Y_G]$.  We can now state our main theorem.

\begin{theorem}\label{th:MainTheoremA}% 
$\Graphs = \Mot$.
\end{theorem}   

The theorem immediately implies that Conjecture~\ref{conjK} is false.
For, if the conjecture were true, all functions of the form $[X]$
would be in $\ring$.  In particular, they would be rational functions.
However, if we let $X$ be the closed subscheme of $\A^1_{\Z}$ defined
by $px=0$ for $p$ a given prime, then $[X](q)=q$ if $p|q$ and $0$
otherwise.  Thus $[X]$ can not be a rational function.  Of course,
other more interesting examples of $X$ such that $[X]$ is not rational
exist.

\subsection{Stanley's Reformulation of Conjecture~\ref{conjK}}
The proof of Theorem~\ref{th:MainTheoremA} is based on Stanley's
reformulation of Kontsevich's conjecture in terms of a polynomial
$Q_G$ which is, roughly speaking, dual to $P_G$.  In \cite{stanley},
Stanley sets
\begin{equation}\label{eq:StanleyPolynomial}
Q_{G } = \sum_{T} \prod_{e\in T} x_e.
\end{equation}
where the sum again runs through all spanning trees.  For $G$
connected, $Q_G$ is homogeneous of degree $\num{E(G)} - b_1(G)$.  Let
$X_G = \A^E - V(Q_G)$.  Stanley showed that Kontsevich's conjecture is
equivalent to the following analogous conjecture:
\begin{conjecture}\label{conjS}
For all graphs $G$, $[X_G]$ is a polynomial.
\end{conjecture}
In fact, we will see in Theorem~\ref{th:graphs} that the
$\ring$-submodule of $\Mot$ generated by the $[X_G]$ is exactly the
same as the one generated by the $[Y_G]$.

The schemes $X_G$ are, however, more tractable than the $Y_G$ ---
particularly when the graph $G$ is simple (i.e., has neither loops nor
multiple edges) and has an apex.  This is because, when $G$ is simple,
the polynomial $Q_G$ has a simple expression as a determinant via
the Matrix-Tree theorem (see section~\ref{matrixtree}).  This
expression simplifies even further when $G$ has an apex.  On the other
hand, while $P_G$ can also be expressed as a determinant, this
expression is combinatorially complicated.

A vertex $v$ is said to be an apex if there is an edge from $v$ to
every other vertex in $G$.  Suppose that $G$ is an arbitrary simple
graph with vertex set $V=\{v_1,\ldots, v_n\}$.  Then we form a graph
$G^*$ with apex by simply adding a vertex $v_0$ and connecting it by
an edge to all other vertices.  All graphs with apex can be obtained
through this process.

Using the Matrix-Tree theorem, Stanley showed that, for any field $K$,
$X_{G*}(K)$ is isomorphic to the set of $n\times n$ nondegenerate,
symmetric matrices $M$ satisfying the condition that
\begin{equation}\label{eq:ZoCondition}
M_{ij}=0\ \text{if $i\neq j$ and there is no edge from $v_i$ to
$v_j$}.
\end{equation}  
Here $i,j\in[1,n]$.

We then let $Z^o_G$ be the scheme of all $n\times n$ nondegenerate,
symmetric matrices $M$ satisfying condition~\ref{eq:ZoCondition}.
(See section~\ref{matrixtree}.)  Stanley's observation essentially
shows that $Z^o_G\cong X_{G^*}$.  Thus, the following conjecture,
stated by Stembridge as Conjecture 7.1 \cite{stembridge}, would follow
from Conjecture~\ref{conjS}.
\begin{conjecture}\label{conjSt}
For every simple graph $G$, $[Z^o_G]$ is a polynomial.
\end{conjecture}

Note that, while Conjectures~\ref{conjK} and~\ref{conjS} are trivial
when $G$ is disconnected, Conjecture~\ref{conjSt} is not.  This is
related to the fact that the operation $G\mapsto G^*$ always produces
a connected graph.

However, we will see that Conjecture~\ref{conjSt} is also false.

For any subgraph $H$ of $G$, let $G-H$ be the graph obtained by
removing the edges in $H$ but leaving all vertices.  Note that
$(G-H)^* = G^* - H$.  If $G$ is a simple graph with $n$ vertices, then
$G$ is contained in the complete graph $K_n$.  We define the {\it
complement} $G^o$ of $G$ to be the graph $K_n -G$.  Note that $(G^o)^*
= (DG)^o$ where $D$ is the operation of adding a disjoint vertex.

It becomes convenient at this point to shift attention from $G$ to its
complement.  We therefore define $Z_G=Z^o_{G^o}$.  When $G$ has
vertices $\{v_1,\ldots, v_n\}$ as above, $Z_G$ is then the scheme of
all $n\times n$ matrices $M$ satisfying the condition
\begin{equation}\label{eq:ZCondition}
M_{ij}=0\ \text{if there is an edge from $v_i$ to $v_j$}.
\end{equation} 

As partial justification that the schemes $Z_G$ are more natural than
the scheme $Z_G^o$, we mention that many of the results obtained thus
far on Conjecture~\ref{conjS} are most easily stated in terms of the
$[Z_G]$.  For example, in Theorem 5.4 of \cite{stanley}, Stanley
showed that Conjecture~\ref{conjS} holds when $G= K_n - K_{1,s}$ where
$K_{1,s}$ is a star (one vertex connected by edges to $s$ other
vertices) and $s\leq n-2$.  In the case $n=s+2$, $G=\Gamma^*$ with
$\Gamma=K_{s+1} - K_{1,s}$.  Thus $\Gamma=K_{1,s}^o$, and $X_G=
Z_{\Gamma}^o = Z_{K_{1,s}}$.  It follows that Stanley's Theorem 5.4 is
equivalent to the statement that $[Z_{K_{1,s}}]\in\Z[q]$.

\subsection{Overview}
Let $\Graphs_*$ be the $\ring$-module generated by all functions of
the form $[Z_G]$ for $G$ a simple graph.  Since
$[Z_G]=[Z_{G^o}^o]=[X_{(G^o)^*}]$, it is clear that
$\Graphs_*\subset\Graphs$.  Therefore the following theorem implies
Theorem~\ref{th:MainTheoremA}:
\begin{theorem}\label{th:MainTheoremB}
$\Graphs_* = \Mot.$
\end{theorem}

The proof of Theorem~\ref{th:MainTheoremB} involves two steps.  In the
first, we study certain incidence schemes $A_G(s,r,k)$.  These schemes
are defined so that, when $K$ is a field, the $K$ points of
$A_G(s,r,k)$ are the set of pairs $(Q,f)$ with $Q$ a symmetric
bilinear form on $K^s$ of rank $r$ and $f$ a function from $V(G)$ to
$K^s$ whose span is of dimension $k$.  The pair $(Q,f)$ is also
subject to the incidence condition that
\begin{equation}\label{ACondition}
Q(f(v_i), f(v_j)) = 0\text{ if there is an edge from $v_i$ to $v_j$.}
\end{equation}

If $G$ has $n$ vertices, then $[A_G(n,n,n)]=[Z_G][\GL_n]$.  Since
$[\Gl_n]\in\ring$, this implies that $[A_G(n,n,n)]\in\Graphs_*$.
Moreover, there are important relations between the $A_G(s,r,k)$ for
varying $s,r$ and $k$, and between the $A_G(s,r,k)$ for varying $G$.
By exploiting these relations, we will see that the $\ring$-module
generated by the $[A_G(s,r,k)]$ is exactly $\Graphs_*$.

This fact allows us to shift our focus from the symmetric form $Q$ to
the function $f$.  In particular, for each $s$ we consider the scheme,
$J_G(s)=\union_k A_G(s,s,k)$.  Again, it turns out that the
$\ring$-module generated by the $J_G(s)$ is exactly $\Graphs_*$.  And
the $J_G(s)$ turn out to be quite manageable schemes because the
dimension of the span of $f$ is allowed to vary.

The second step in our proof of Theorem~\ref{th:MainTheoremB} involves
comparing the $J_G(s)$ to the representation spaces of matroids.  For
any matroid $M$, we define a scheme $X(M,s)$.  For $K$ a field,
$X(M,s)(K)$ is the set of all possible representations of $M$ in
$K^s$.  We then let $\Matroids$ denote the $\ring$-module generated by
all functions $[X(M,s)]$.  As we will see in Section~\ref{sec:mnev},
it essentially follows from Mn\"ev's Universality
Theorem~\cite{Mnev} that $\Matroids = \Mot$.  On the other hand, we
prove that, for each matroid $M$, there is a finite set of graphs
$\{G_i\}$ and rational functions $a_i\in \ring$ such that
\begin{equation}\label{eq:IntroPunchline}%
[X(M,s)] = \sum a_i [J_G(s)].
\end{equation} 

This equation proves that $\Matroids\subset\Graphs$ and, thus, it
proves Theorem~\ref{th:MainTheoremB}.  Moreover, as we will see,
(\ref{eq:IntroPunchline}) can be used even without Mn\"ev Universality
to produce a contradiction to Conjecture~\ref{conjSt}.  This is
because there are matroids $M$, for example the Fano matroid, which
are representable only over fields of characteristic $2$.  Thus, for
such matroids, $[X(M,r)]$ (with $r$ equal to the rank of $M$) could
not possibly be a rational function as Conjecture~\ref{conjSt} and
(\ref{eq:IntroPunchline}) would demand.  As Conjecture~\ref{conjK}
implies Conjecture~\ref{conjSt}, this shows that
Conjecture~\ref{conjK} is false.

\subsection{Forest Complements}

A considerable amount of work has been done to find examples of graphs
for which Conjecture~\ref{conjK} (resp. Conjecture~\ref{conjS},
Conjecture~\ref{conjSt}) holds and to compute the functions $[Y_G]$
(resp. $[X_G]$, $[Z_G^o]$) explicitly \cite{chung, stanley,
stembridge, yang}.  It remains an interesting question to determine
the largest classes of graphs for which these conjectures remain true.

The class of graphs for which Conjecture~\ref{conjS} is true is
already known to include to include various interesting graphs.
Stanley showed that $X_{K_n - K_m}$ is polynomially countable.  Chung
and Yang then computed the polynomial $[X_{K_n - K_m}]$ explicitly
\cite{chung}.  Yang showed that $X_G$ is polynomially countable when
$G$ is an outplanar graph.  And, as mentioned above, Theorem 5.4 of
\cite{stanley} shows that $Z_{K_{1,s}}$ is polynomially countable.  In
fact, the theorem is equivalent to the statement that $Z_G$ is
polynomially countable when $G$ is the union of a star and a discrete
graph.  (Stanley also computes $[Z_G]$ explicitly in this case.)

Recall that a {\em forest} is a graph with no cycles.  
In section~\ref{forests}, we show that $Z_F$ is polynomially countable
whenever $F$ is a forest.  This
generalizes Stanley's Theorem 5.4 and implies
that Conjecture~\ref{conjSt} holds for forest complements.  
The result is essentially a
consequence of the manageability of the schemes $J_F(s)$ which allows
us to compute $[J_F(s)]$ inductively in terms of the $[J_{F'}(s)]$
for smaller forests $F'$.  

\subsection{Acknowledgments}   We must express our deep appreciation
to T.~Chow and N.~Fakhruddin for reading several versions of the early
stages of this work.  Without their support, it is unlikely that we
would have had the tenacity to push this work to its final conclusion.
We also owe a deep debt of gratitude to B.~Totaro who pointed us to
Mn\"ev's Universality Theorem.  In addition we would like to thank
F.~Chung, W.~Fulton, D.~Kreimer, N.~Mn\"ev, M.~Nori, R.~Stanley,
J.~Stembridge, B.~Sturmfels and C.~Yang for useful conversations and
email correspondences.

The present paper is an updated version of our Max-Planck-Institut
f\"ur Mathematik preprint~\cite{bb}.  The second author would like to
thank the Max-Planck-Institut f\"ur Mathematik in Bonn for providing
the wonderful environment in which most of his work on \cite{bb} was
done.  He would also like to thank his fellow visitors to Bonn,
especially J.~Furdyna, R.~Joshua, A.~Knutson, S.~Lekaus, P.~Mezo,
A.~Schwarz, D.~Stanley and B.~Toen, for useful and encouraging
conversations.

\section[graphs]{Preliminary Results}

In this section we carry out two minor adjustments to two theorems
of Stanley.  

\subsection{The Module of all Graphs}  
The first adjustment is an amplification of Proposition 2.1 of
\cite{stanley}.  It concerns the relation between the schemes $Y_G$
and the schemes $X_G$.

\begin{proposition}\label{pr:graphs} The subgroup of 
$\Mot^+$ generated by the functions $[X_G]$ is equal to 
the subgroup generated by the functions $[Y_G]$.
\end{proposition}

We remark that the proof of this proposition is completely contained in 
Stanley's proof of his Proposition 2.1.  However, for the convenience of
the reader, we translate Stanley's proof into our own setting.

\begin{proof}
Let $S$ be a subset of $E=E(G)$.  Let $\A^S$ be the image of the obvious
inclusion of $i^S:\A^{\num{S}}\to A^E$.  Let $\G_m^S=i^S(\G_m^{\num{S}})$.
Note that, as $S$ varies over all subsets of $E$, the subschemes
$\G_m^S$ stratify $\A^E$.

For any subscheme $X\subset\A^E$, let $X_S=X\cap\A^{E-S}$ (resp. 
$X_S^+=X\cap\G_m^{E-S})$.  Thus $X_S$ is the intersection of $X$ with 
the hyperplanes defined by the equations $x_e$ for $e\in S$.  
Note that $X_{\emptyset}=X$, and, as $S$ varies over the subsets of $E$, 
the subschemes $X_S^+$ stratify $S$.  
We therefore have, 
\begin{equation}\label{eq:Decomp}
[X_S]=\sum_{T\supset S} [X_T^+]
\end{equation}
and, by the Inclusion-Exclusion Principle,
\begin{equation}\label{eq:IncEx}
[X_S^+]=\sum_{T\supset S} (-1)^{\num{(T-S)}} [X_T].
\end{equation}

By inspecting the $Q_G$, it is easy to see that $X_{G,S}\cong X_{G-S}$
and $X_{G,S}^+\cong X_{G-S}^+$.  Dually, if $S$ is a forest,
$Y_{G,S}\cong Y_{G/S}$ (resp. $Y_{G,S}^+\cong Y_{G/S}^+$) where $G/S$
is the graph obtained by contracting each component of $S$ to a point.
On the other hand, if $S$ is not a forest, it is easy to see that
$Y_{G,S}$ is empty.

Now, as Stanley notes, 
$Q_G(x) = P_G(1/x)\prod_{e\in E} x_e$.  Thus 
\begin{equation}\label{eq:XplusY}
X_{G,\emptyset}^+\cong Y_{G,\emptyset}^+
\end{equation} through the map 
$x\mapsto 1/x$.  

Putting our equations together we obtain the following:
\begin{eqnarray}
\lb Y_G \rb &=& \sum_{S\subset E\atop b_1(S)=0}
                \sum_{T\subset G/S} (-1)^{\num{T}} 
                \lb X_{(G/S)-T} \rb,\\
\lb X_G \rb &=& \sum_{S\subset E}\sum_{T\subset E-S\atop b_1(T)=0}
                (-1)^{\num{T}} \lb Y_{(G-S)/T} \rb. 
\end{eqnarray}
Together, these two equations, the first of which appears (in a different
notation) as Proposition 4.1 of \cite{stembridge}, prove the proposition.
\end{proof}

The proposition implies the following theorem as a corollary.
\begin{theorem}\label{th:graphs} $\Graphs$ is equal to the 
$\ring$-submodule of $\Mot$ spanned by the $[X_G]$.
\end{theorem}

We remark that, as $[Y_G]=q^n-q^{n-1}$ for $G$ a cycle of length $n$,
$\ring$ is itself a submodule of $\Graphs$.

\subsection{An Observation on Polynomial Countability}  Our second 
adjustment to Stanley's results is to Proposition 2.2 of
\cite{stanley}.  This proposition, which Stanley deduces from the Weil
conjectures, essentially states that, if $X$ is a scheme of finite type
over $\Z$, then the knowldege that $[X]\in\Q[q]$ implies that, in fact,
$[X]\in\Z[q]$.  

In Section~\ref{forests}, we require a result which is analogous to
Stanley's Proposition 2.2 but easier to prove.  While the result is
not strictly weaker than Stanley's proposition, it does not require
the Weil conjectures.  Rather, it is a consequence of the Euclidean
algorithm and the infinitude of the primes.

\begin{proposition}\label{th-prational} 
If $f\in \ring$ and $f(q)\in\Z$ for all $q\in\primaries$, then
$f\in\Z[q]$.
\end{proposition}

We will use the proposition in the case $f=[X]$ for 
$X$ a scheme of finite type over $\Z$.

\begin{proof} Write $f = a/s$ with $a\in\Z[q]$ and $s\in\calS$.  
Since $s$ is monic, we can write $f = d + r/s$ with
$d,r\in\Z[q]$ and $\deg(r)<\deg(s)$.  But this implies that 
$r(q)/s(q)\in\Z$ for all $q$ which implies that $r=0$.   Thus $f=d$.
\end{proof}

We remark that it is possible, using the rationality of 
the Zeta function, proved by Dwork in \cite{dwork},
to prove a much stronger result which implies both Stanley's 
Proposition 2.2 and our Proposition~\ref{th-prational} above.

\begin{theorem} Let $p$ be a prime and $X$ and $Y$ schemes of finite type
over $\F_p$.  For all $n>0$, let $f_X(p^n) = \num{X(\F_{p^n})}$ 
(resp. $f_Y(p^n) = \num{Y(\F_{p^n})}$).
Let $g\in\Q(q)$ be a rational function.
\begin{enumerate}
\item[(a)] 
If $f_X(p^n)=g(p^n)$ for almost all $n$, then $g\in\Z[q]$,
and $f_X(p^n)=g(p^n)$ for all $n$.
\item[(b)] If $f_X(p^n)=f_Y(p^n)$ for almost all $n$, then 
$f_X(p^n)=f_Y(p^n)$ for all $n$.
\end{enumerate}
\end{theorem}

We do not need this result because the denominators we
consider always lie in the multiplicative set $\calS$.  The theorem
would be useful, however, if this were not the case.  The proof
is left to the interested reader.
 
\section[schemes]{Determinantal Schemes}\label{schemes}

In this section, we collect certain basic properties of determinantal
schemes which are necessary for the rigorous definition of the
incidence schemes $A_G(s,r,k)$ discussed in the introduction and for
some of the discussion in Section~\ref{sec:mnev}.  We first describe
the general theory of determinantal schemes in functorial language and
then restrict to the specific case of determinantal schemes over $\Z$
that are the focus of the paper.

The results of this section are not strictly necessary for the proof
that Conjecture~\ref{conjK} is false.  In fact, the proof of
Theorem~\ref{th:MainTheoremB} does not require the fact that the
$A_G(s,r,k)$ are schemes.  It is enough to consider the $A_G(s,r,k)$
as functions from finite fields to sets, and for this the definition
given in the introduction suffices.  Therefore, the reader only
interested in the proof that Conjecture~\ref{conjK} is false can
safely skim this section.

\subsection{General Theory} 
Let $S$ be a scheme and let $E$ and $F$ be locally free
$\clO_S$-modules of ranks $e$ and $f$ respectively.  Write
$\Hom_{\clO_S}(E, F)$ for the abelian group of homomorphisms from $E$ to
$F$.  The scheme of homomorphisms $\Hom(E,F)$ is then an abelian group
scheme over $S$ representing the functor
\begin{equation}
T\cmapsto Hom_{\clO_T}(E_T, F_T).
\end{equation}

Now let $r\in\N$, and write $\Hom_{\clO_S, r}(E,F)$ for the set of
morphisms $$\phi\in\Hom_{\clO_S}(E,F)$$ such that $\coker\phi$ is
locally free of rank $f-r$.  The functor
\begin{equation}
T\cmapsto Hom_{\clO_T,r}(E_T, F_T).
\end{equation}
is representable by a scheme $\Hom_{r}(E,F)$.  

To see this, first note that we can assume without loss of generality
that the base $S$ is affine and the sheaves $E$ and $F$ are both free.
Let $S=\Spec A$ and let $T=\Spec B$.  Then $\Hom_{\clO_T,r}(E_T,F_T)$
is equal to the set of maps $\phi:B^e\to B^f$ such that $\coker\phi$
is a projective $B$-module of rank $f-r$.  Now, $\coker\phi$ is
projective of rank $f-r$ if and only if for every prime $\wp\in\Spec B$
the rank of $\phi\otimes (B_{\wp}/\wp)$ is $r$.  And this will be the
case if and only if every $(r+1)\times (r+1)$ minor in $\phi$
vanishes, but some $r\times r$ minor is invertible.

Let $\displaystyle{\{y_{ij}\}_{i=1}^e}_{j=1}^f$ 
be a set of formal variables, and consider each $y_{ij}$ 
as an entry in an $e\times f$ matrix.  Let $A[y]$ be the polynomial
ring in all variables $y_{ij}$.  
For each $k$, let $m_l^k\in\Z[y_{ij}]$
be a complete list of the $k\times k$ minors, and let $I_k$ be the 
ideal generated by the $m_l^k$. 
In this notation, 
$Hom_r(\clO_S^e, B^f)$ is the locally closed subscheme of 
$Hom(\clO_S^e, \clO_S^f)$
given by the union of the affine schemes
\begin{equation}\label{eq:AffineHomr}
\union_i \Spec (A[y]/I_{r+1})_{(m_i)}.
\end{equation}
It follows that the set of points associated to $Hom_r(\clO_S^e, \clO_S^f)$ 
is simply $\displaystyle\cap_{i} V(m_i^{r+1}) - \cap_{i}V(m_i^r)$.  

Thus the schemes $Hom_r(E,F)$ form a stratification of $Hom_r(E,F)$.  
To see this stratification in a coordinate free manner, we consider
the {\em determinantal sets} $Hom_{\clO_S,\leq r}(E,F)$.  This is the set 
of all maps $\phi\in\Hom_{\clO_S}(E,F)$ such that the stalk 
of $\coker\phi$ has 
rank less than or equal to $r$ at every point $x\in S$.   The
functor
$$
T\cmapsto\Hom_{\clO_T, \leq r}(E_T,F_T)
$$
is then represented by the {\em determinantal scheme} 
$\Hom_{\leq r}(E,F)$.  These scheme have been studied extensively. 
\cite{HEagon, Room}. 

In the case that $S=\Spec A$, $E=A^e$ and $F=A^f$ considered above, 
$\Hom_{\leq r}(E, F)$ is simply $\Spec A/I_{r+1}$.  Thus 
\begin{equation}
\Hom_{r}(E,F) = \Hom_{\leq r} (E,F) - \Hom_{\leq r-1}(E,F).
\end{equation}

\subsubsection{Maps to the Grassmanian}
Write $\Gr(r,E)$ for the Grassmanian of $r$ planes in $E$.  
$\Hom_r(E,F)$ is that this scheme is equipped with two maps to
Grassmanians.  We have a map $p:\Hom_r(E,F)\to\Gr(r, F)$ given
essentially by sending $\phi$ to its image.  And we have a map
$q:\Hom_r(E,F)\to\Gr(e-r, F)$ given by sending $\phi$ to the its
kernel.

\subsubsection{Function Spaces} When $V$ is a finite set we write 
$\Fun (V,E)$ for $\Hom (\clO_S^V, E)$ (resp. $\Fun_r(V,E))$ for
$\Hom_r(\clO_S^V, E)$.

\subsubsection{Symmetric Bilinear Forms}  Let $E^{\vee}$ denote
the dual of $E$.  There is a natural transpose automorphism
\begin{equation}
t: \Hom(E,E^{\vee})\to\Hom(E,E^{\vee})
\end{equation}
and we define $\Sym E$ to be the subscheme fixed by $t$.  
We then write $\Sym_r E$ (resp. $\Sym_{\leq r} E$) 
for the scheme-theoretic intersection
of $\Sym E$ with $\Hom_r(E,E^{\vee})$ (resp. $\Hom_{\leq r}(E,E^{\vee})$.

\subsection{A Specific Case}  We will be primarily interested in the 
case $S=\Spec\Z$, $E=\clO_S^e$ and $F=\clO_S^f$.  In this case,
$\Hom(E,F)$ is $\Spec\Z[y]$.  $\Hom_{\leq r}(E,F)$ is the closed
subscheme in $\Hom(E,F)$ defined by the $(r+1)\times (r+1)$ minors.
And $\Hom_r(E,F)$ is the Zariski open subset of $\Hom_{\leq r}(E,F)$
defined by requiring at least one $r\times r$ minor to be invertible.

These equalities can be used without reference to the preceding
general theory to define the $\Hom_r(E,F)$.  It follows directly that
for any field $K$, $\Hom_r(E,F)(K)$ is the set of maps from $K^e$ to
$K^f$ of rank $r$.

Similarly, when a $E=\clO^e_S$ with $S=\Spec\Z$, $\Sym E$ can be
viewed as the closed subscheme of $\Z[y]$ defined by the equations
$y_{ij}=y_{ji}$.  $\Sym_{\leq r} E$ is then the closed subscheme of
$\Sym E$ defined by the $(r+1)\times (r+1)$ minors.  And $\Sym_r E$ is
the Zariski open subset of $\Sym_{\leq r} E$ defined by requiring at
least one $r\times r$ minor to be invertible.  The $K$ points of
$\Sym_r E$ are the bilinear forms on $K^e$ of rank $r$.

\subsection{Polynomial Countability}
With $E$ trivialized over $\Spec\Z$ as above, we write $\Gl_e$ for
$\Hom_e(E,E)$, $\Gr(r, e)$ for $\Gr(r, E)$, and $\Sym_r^e$ for 
$\Sym_r E$.

We now
list a few results concerning the polynomial countability of the
schemes just discussed.

\begin{eqnarray}
\lb\Gl_n\rb       &=& (q^n -1)(q^n - q) \cdots  (q^n - q^{n-1})\\
\lb\Gr(a, b)\rb &=& \frac{\lb\Gl_{b }\rb }{\lb\Gl_a\rb\lb\Gl_{b-a}\rb 
                      q^{a(b-a) }  }\\ 
\lb\Hom_r(E,F)\rb &=& \lb\Gr(r, e)\rb\lb\Gr(r,f)\rb\lb\Gl_r\rb
\end{eqnarray}

The first two of the above equalities are well known and the last is easy.  
Note that each of the functions given is a polynomial lying in 
the multiplicative set $\calS$.

A more difficult formula is the following, taken from 
\cite{macwilliams}.
\begin{equation}\label{eq:mwa}
\lb\Sym^n_r \rb  = \left\{\begin{array}{ll}
           \prod_{i=1}^s\frac{q^{2i}}{q^{2i}-1}\cdot
           \prod_{i=0}^{2s-1}(q^{n-i} -1),\ 0\leq r=2s\leq n,\\
           \prod_{i=1}^s\frac{q^{2i}}{q^{2i}-1}\cdot
           \prod_{i=0}^{2s}(q^{n-i} -1),\ 0\leq r=2s+1\leq n 
                           \end{array}\right.
\end{equation}

Note again that $[\Sym^n_r]\in\calS$.

\section[matrix tree]{The Matrix Tree Theorem}\label{matrixtree}

Stanley's positive results mentioned in the introduction were mainly
consequences of the {\it Matrix-Tree Theorem} of Kirchhoff, Borchardt
and Sylvester which gives an expression of the polynomial $Q_G$ as the
determinant of a symmetric matrix.  As this theorem is also basic to
our results, we describe it in this section after fixing some useful
notation.

\subsection{Notation}\label{notation}
When $G$ is a simple graph, an assumption we will make
for the remainder of this paper, $E$ can be considered as a subset
of $\Sym^2 V$.  For $v,w\in V$, we write $e_{vw}$ for the set 
$\{v,w\}$.  Thus the statement $e_{vw}\in E$ means that there is an
edge in $G$ connecting $v$ to $w$.   

It is convenient to pick an ordering $V=\{v_1,\ldots, v_{n_G}\}$ of the
edges $V$.  Thus $n_G=\num{V(G)}$.  We write $n$ for $n_G$
when there is only one graph under consideration.

Set $e_{ij}=e_{v_iv_j}$.  We write $x_{ij}$ for the variable
$x_{e_{ij}}$ when $e_{ij}\in E$, and we extend this notation by
setting $x_{ij}=0$ when $e_{ij}\not\in E$.

\subsection{The Laplacian}\label{laplacian}
Let $L = L_{ij}$ be the $n\times n$ matrix defined
by
$$
L_{ij} = \left\{\begin{array}{r@{\quad}l}
        \sum_{k=1}^n x_{ik} & \hbox{if $i=j$}\\          
                    -x_{ij} & \hbox{if $i\neq j$}      
                 \end{array}\right.
$$
Let $L_0$ be $L$ with first row and the first column removed.
$L$ is called the {\it generic Laplacian matrix of $G$} and $L_0$
the {\it reduced generic Laplacian}.  The following theorem 
can be found in the work of Cayley, Kirchhoff,  Maxwell and Sylvester.  
For a proof, see \cite{StanleyEnComb2}.

\begin{theorem}[The Matrix-Tree Theorem]\label{matrix-tree} 
$Q_G = \det L_0$.
\end{theorem}

Now, as in the introduction, let $Z_G^0$ be the scheme of all 
$n\times n$ symmetric, non-degenerate bilinear forms $M_{ij}$ such
that $M_{ij}=0$ whenever $i\neq j$ and $e_{ij}\not\in E$.  In the notation
of section~\ref{schemes}, $Z_G^0$ is simply the closed subscheme of 
$\Sym^n_n$ defined by the equations $y_{ij}=0$ for all $i\neq j$ with
$e_{ij}\not\in E$.   

Our use of Theorem~\ref{matrix-tree}, is based on the following important
consequence, recognized by Stanley.

\begin{theorem}\label{th-stanley} $X_{G^*}\cong Z_G^o$.
\end{theorem}
\begin{proof}  Let $\Z[x]$ be the ring generated by the variables 
$x_{ij}$ for $0\leq i<j\leq n$.  Let $I$ be the ideal generated by
the variables $x_{ij}$ for all pairs $i<j$ with $e_{ij}\not\in E$.
Then $X_{G^*}=\Spec A$ with $A=(\Z[x]/I)_{Q_G}$.

On the other hand, let $\Z[y]$ is the ring generated by all $y_{ij}$ for
$i,j\in\{1,\ldots, n\}$, and let $J$ be the ideal generated by all
expressions of the form $y_{ij}-y_{ji}$ for $i\neq j$ and $y_{ij}$ for
$i\neq j$ and $e_{ij}\not\in E$.  Then, letting $D$ be the determinant
of the matrix of $y_{ij}$'s, $Z_G^0=\Spec B$ with 
$B=(\Z[y]/J)_D$. 

Let $p:\Z[y]\to\Z[x]$ be the map 
\begin{equation}
y_{ij}\mapsto\left\{\begin{array}{r@{\quad}l}
                   \sum_{k<i} x_{ki} + \sum_{i<k} x_{ik} & i=j\\
                   -x_{ij}                               & i< j\\
                   -x_{ji}                               & j<i 
                  \end{array}\right.
\end{equation}

Let $q:\Z[x]\to \Z[y]$ be the map 
\begin{equation}
x_{ij}\mapsto\left\{\begin{array}{r@{\quad}l}
                  \sum_k y_{jk}                          & i=0\\
                  -y_{ij}                                & i>0
                    \end{array}\right.
\end{equation}

It is easy to verify that $p(I)\subset J$, that $q(J)\subset I$, 
and that $p$ and $q$ give inverse isomorphisms between the rings
$\Z[x]/I$ and $\Z[y]/J$.  It then follows from the Matrix-Tree theorem
that $p(Q_G)=D$.  Thus $p$ and $q$ give inverse isomorphisms between 
the rings $A$ and $B$.
\end{proof}

As mentioned in the introduction, it is convenient to shift our attention
from the simple graph $G$ to its complement.  We therefore set 
$Z_G= Z_{G^o}^o$.  Thus $Z_G$ is the subscheme of $\Sym^n_n$ defined
by the equations $y_{ij}=0$ for every pair $i,j$ with $e_{ij}\in E$.
And $Z_G=X_{(G^o)^*}=X_{(DG)^o}$ where $D$ is the operation of adding
a disjoint vertex.

\begin{example} Let $G$ be a graph with $n$ vertices and
no edges.  Then $Z_G\cong\Sym^n_n$.  This is recognized in
\cite{stanley}.  By Equation \ref{eq:mwa}, it follows that
$[Z_G]\in\Z[q]$.  In fact, $[Z_G]\in\calS$, and this shows that
$\ring\subset\Graphs_*$.
\end{example}

\section[incidence schemes]{Incidence Schemes}\label{sec-incidenceschemes}

We now introduce the incidence schemes mentioned in the introduction.
At first, we work in full generality over a base scheme $S$.  But
our main interest is the case $S=\Spec\Z$.

\begin{definition} Let $W$ be a locally free $\clO_S$-module, and let
$G$ be a graph.   We write $A_G(W)$ for the closed subscheme of 
$\Sym W\times\Fun(V,W)$ consisting of pairs $(Q,f)$ satisfying the
condition that 
\begin{equation}
Q(f(v),f(w))=0\ \text{if $e_{vw}\in E(G)$.}
\end{equation}
If $r$ and $k$ are integers, we write $A_G(W,r,k)$ for 
$$
A_G(W)\cap\left(\Sym_r(W)\times\Fun_k(V,W)\right).
$$
When $S=\Spec\Z$ and $W=\clO_S^s$, we write $A_G(s)$
for $A_G(W)$ and  
$A_G(s,r,k)$ for $A_G(W,r,k)$.
\end{definition}

The $A_G(W,r,k)$  form a stratification
of $A_G(W)$  by locally closed subschemes.
Note that $A_G(s,r,k)$ is empty unless 
$0\leq k <n$ and $0\leq r,k \leq s$.
Also note that $A_G(s,r,0) = \Sym^s_r$, and 
$A_G(W,0,k)=\Fun_k(V,W)$.  Thus  
$[A_G(s,r,0)]$ and $[A_G(s,0,k)]$ are both in $\Z[q]$.

Now assume that $V(G)=\{v_1,\ldots, v_n\}$ as in paragraph~\ref{notation}.
Recall from the introduction that $\Graphs_*$ is the $\ring$-submodule
of $\Mot$ spanned by the functions $Z_G$.

\begin{theorem}\label{th-KGmandXG}
\begin{enumerate}
\item[(a)] $A_G(n,n,n) \cong Z_G \times \Gl_{n}$.  
\item[(b)] 
$\Graphs_*$ is exactly equal to the $\ring$-module
generated by the functions $[A_G(n,n,n)]$.
\end{enumerate} 
\end{theorem}
\begin{proof} We first remark that (b) follows directly 
from (a) and the fact that $[\Gl_n]\in\calS$.

To prove (a) we let $W=\clO^n_S$ with $S=\Spec\Z$.  Then 
$\Fun_r(V,W)=\Gl_n$.  The map $(Q,f)\mapsto (f^tQf, f)$ then
identifies $A_G(n,n,n)$ with $Z_G\times\Gl_n$.  (Here $f^t$ denotes
the transpose of $f$).
\end{proof}

\begin{remark} Let $Z_G(r)$ be the scheme consisting of all $n\times n$
symmetric bilinear forms of rank $r$ such that $M_{ij}=0$ whenever
$e_{ij}\in E$.  These schemes have been studied implicitly in
\cite{chung, stanley}.  In Stanley's notation, $h(G,r)=[Z_{G^o}(r)]$,
and Chung and Yang call a graph $G$ {\em strongly admissible} if
$Z_{G^o}(r)$ is polynomially countable for all $r$.  An easy
modification of the proof above shows that $A_G(n,r,n)\cong
Z_G(r)\times\Gl_n$.
\end{remark}

\section[macwilliams]{Extensions of Bilinear Forms}\label{sec-macwilliam}

In this section, we review a result of MacWilliams~\cite{macwilliams}
counting the number of ways to extend a bilinear form of rank $r_1$ to 
a bilinear form of rank $r_2$.  This count will be important in the next
section for finding relations among the $A_G(s,r,k)$.

Let $Q$ be a fixed bilinear form on $\F_q^{d_1}$ with rank $r_1$.  Let
$C_Q(d_2, r_2, d_1 , r_1)$ be the number of ways to extend $Q$ to a
form on $\F_q^{d_2}$ of rank $r_2$.  The following result is Lemma 4
of \cite{macwilliams}.
\begin{theorem}\label{th-macwilliams} 
$$
C_Q(d_1+1, r_2 ,d_1, r_1)=\left\{\begin{array}{r@{\quad}l}
                        q^{r_1} & r_2=r_1       \\
            q^{r_1+1} - q^{r_1} & r_2=r_1+1     \\
        q^{d_1 + 1} - q^{r_1+1} & r_2 = r_1+2   \\
                            0 & \hbox{otherwise.}
                       \end{array}\right.
$$
\end{theorem}

Note that $C_Q(d_1+1, r_2 ,d_1, r_1)$ only depends on $d_1, r_2$ and
$r_1$.  By induction on $d_2-d_1$, we can show that $C_Q(d_2, r_2,
d_1, r_1)$ only depends on the integer parameters $d_2, r_2, d_1$ and
$r_1$.  Thus we simply write $C_(d_2, r_2, d_1, r_1)$ for this number.
We can also see by induction that the following recursion is
satisfied
\begin{equation}\label{recursion}
\begin{split}
C(d_2, &r_2, d_1, r_1) = \\
&\sum_{j=0}^2 C(d_2,r_2, d_1+1, r_1+j)
C(d_1+1,r_1 +j ,d_1, r_1).
\end{split}
\end{equation}

\begin{corollary}\ 

\begin{enumerate}
\item[(a)] $C(d_2, r_2, d_1, r_1)$  
is a polynomial in $q$. 
\item[(b)] $C(d_2, r_2, d_1, r_1)\neq 0$ iff $d_2\geq r_2$, $d_1\geq r_1,$ 
and $0\leq r_1\leq r_2\leq r_1 + 2(d_2 - d_1)$.
\end{enumerate}
\end{corollary}
\begin{proof} (a) follows directly from the recursion formula~\ref{recursion}.

The necessity of the first two inequalities of (b) are obvious for
dimension reasons (the rank of a bilinear form can not be greater than
the dimension of the ambient space).  Necessity of the third
inequality follows from formula~\ref{recursion} by induction.

We prove the sufficiency of the the inequalities in (b) by
induction on $i=d_2 - d_1$ using formula~\ref{recursion}.
We do not actually need this for the rest of the paper so the 
reader may safely skip the proof.

For $i=0$ sufficiency is obvious.  For $i=1$ the sufficiency results
from the fact that $C(d_1+1, r_2, d_1, r_1)\neq 0$ 
iff $r_2\in [r_1, r_1 +2]$
when $d_1 \neq r_1$ and iff $r_2\in [r_1, r_1 +1]$ when $d_1=r_1$.

Now suppose sufficiency is known
for $d_2 - d_1<i$ and assume that
$(d_2, r_2, d_1, r_1)$ satisfies the conditions in (b)
with $d_2 = d_1 + i$ and $r_2 = r_1 + k$.
By formula~\ref{recursion}, 
$C(d_2, r_2, d_1, r_1)\neq 0$ if there is a $j$ such that both
\begin{enumerate}
\item $C(d_1 + i, r_1+k, d_1 + 1, r_1 +j)\neq 0$ and 
\item $C(d_1 +1, r_1+j, d_1, r_1)\neq 0$.
\end{enumerate}
One computes that (2) is satisfied whenever $j\leq d_1 - r_1 + 1$.
Using the the induction hypothesis, we see that (1) is satisfied
for 
$$
k-2i + 2\leq j \leq \min(k, d_1 - r_1 + 1).
$$
So we need only show that $k-2i + 2\leq \min(k,d_1 -r_1 +1)$.
That $k-2i +2\leq k$ only says that $i\geq 1$ which we are of course
assuming.  And $k-2i + 2\leq d_1 - r_1 +1$ iff
$(d_2 - d_1) + (d_2 - r_2) \geq 1$ which then follows from the 
fact that $d_2\geq r_2$.
\end{proof}

\section[reductions]{Reduction Formulae}\label{sec-reductionformulae}

In this section we give three formulae which allow us to reduce
questions about $A_G(s,r,k)$ for given $s,r$ or $k$ to questions where
$s,r$ or $k$ is smaller.  We also give a formula that allows us to
connect the $A_G(s,r,k)$ to the $A_{DG}(s,r,k)$ where $DG$, as in the
introduction, is the graph obtained from $G$ by adding a disjoint
vertex.

\begin{theorem}\label{firstred}
$$ [A_G(s,r,k)] = [\Gr(k, s)]\sum_j C(s, r, k, j)[A_G(k, j, k)].$$
\end{theorem}
\begin{proof}  
In this proof and the next one we will pick a base field 
$\F_q$ at the beginning and then, for any scheme $X$ we encounter,
write $X$ instead of $X(\F_q)$. 

Write $W=\F_q^s$.  For every map $f\in W^V$ let $\spn{f}$ denote the
span of the $f(v_i)$.  The map $(Q,f)\mapsto \spn{f}$ fibers the set
$A_G(s,r,k)$ over $\Gr(k,s)$.  The fiber over a subspace $U\subset W$
is then the set $A_G(s,r,U)$ of $(Q,f)\in A_G(s,r,k)$ such that
$\spn{f} = U$.  The transitivity of the $\Gl_s$ action on $\Gr(k,s)$
shows that the fibers all have the same number of points.  Thus for any
given $U$
\begin{equation}\label{frone}
\num{A_G(s,r,k)} = \num{\Gr(k,s)}\cdot\num{A_G(s,r,U)}.
\end{equation}

Now let $A_G(s,r,U,j)$ be the set of $(Q,f)\in A_G(s,r,U)$ such
that $Q|_U$ has rank $j$.  This decomposes $A_G(s,r,U)$ into disjoint
subsets.
Consider the map 
\begin{eqnarray*}
p_U: A_G(s,r,U,j) &\to     & A_G(U,j,k)\hbox{ given by}\\
          (Q,f) &\mapsto & (Q|_U, f).
\end{eqnarray*}
The fiber of $p_U$ above a given $(\overline{Q},f)\in A_G(U,j,k)$ is 
$C_{\overline{Q}}(s,r, k, j)$.  Thus
\begin{equation}\label{frtwo}
\num{A_G(s,r,U,j)} = C(s,r,k,j)\cdot\num{A_G(k,j,k)}
\end{equation}
Summing over all the $j$ in Equation~\ref{frtwo} and substituting
the result into Equation~\ref{frone} we obtain the desired result.
\end{proof}

\begin{theorem}\label{secondred}
 $
[A_G(s,r,k)] = $
$$ 
         [Gr(r,s)] \sum_l [\Gr(n-k,n-l)][\Gr(k-l,s-r)][\Gl(k-l)]q^{l(s-r)}[A_G(r,r,l)]
$$
\end{theorem}
\begin{proof}
Write $W=\F_p^s$ and
let $\Psi:A_G(W,r,k)\to \Gr(s-r, W)$ be the
map associating to every $(Q,f)$ the kernel $\ker Q$ of $Q$.  The fiber
of $\Psi$ over a subspace $U\subset  W$ is the set
$A_G(W,r,k)_U$ consisting of all $(Q,f)\in A_G(W,r,k)$ with $Q|_U=0$.  
The transitivity of the action of $\Gl(W)$ on $\Gr(s-r, W)$ 
show then that  
\begin{equation}\label{srone}
\num{A_G(W,r,k)} = \num{Gr(s-r, W)}\cdot\num{A_G(W,r,k)_U}
\end{equation}

Let $T=W/U$
and let $\pi: W\to T$ be the quotient map.  $Q$ reduces in an obvious way
to a form $\overline{Q}$ on $T$.  In fact 
$Q\mapsto \overline{Q}$ is a one-one correspondence
between bilinear forms on $W$ with kernel $U$ and non-degenerate bilinear
forms on $T$.

Now stratify $A_G(W,r,k)_U$ by the dimension of $\spn{\pi\circ f}$. This stratum
corresponding to span $=l$ maps to $A_G(T,r,l)$ by sending $(Q,f)$
to $(\overline{Q}, \pi\circ f)$.     
The fiber above a pair $(\overline{Q},g)$ is identified with set of 

$f: V\to W$
such that $\pi\circ f = g$ and $\spn{f}$ is of dimension $k$. It is elementary
linear algebra to verify that this is given by
$$
[\Gr(n-k,n-l)][\Gr(k-l,s-r)][\Gl(k-l)]q^{l(s-r)}.
$$

Hence putting these fibrations together we get the desired result.
\end{proof}

It is worth recording an important special case of this result.
\begin{corollary}\label{cor-secondred}
$$
[A_G(s,r,s)]=[Gr(n-s,n-r)][Gl(s-r)]q^{r(s-r)}[A_G(r,r,r)].
$$
\end{corollary}
\begin{proof}
To get a non-zero contribution  corresponding to $l$ in the previous theorem need
\begin{enumerate}
\item $r\leq s$.
\item $l \leq r$.
\item $l \leq k$.
\item $l \geq k+r-s$.
\end{enumerate}
In the case of the corollary s=k, so we get $l\leq r$ and $l\geq r$. Hence $l=r$, and the formula reduces to exactly the above.
\end{proof}

We now give a reduction theorem relating the incidence schemes
of $DG$ to those of $G$.

\begin{theorem}\label{Dreduction}
$$
[A_{DG} (s,r,k)] = 
q^k [A_{G}(s,r, k)] + (q^s - q^{k-1}) [A_G(s,r,k-1)].
$$
\end{theorem}
\begin{proof} Let $W=\F_q^s$.
Let $f(V(DG))\to W$ with $\spn{f}$ a $k$-dimensional subspace.
The span of $f|_{V(G)}$ is either a $k$ or a $k-1$ dimensional
subspace.  If $\{v\}=V(DG)-V(G)$, counting the 
possibilities for $f(v)$ proves the lemma.
\end{proof}

\section[motives]{The Module of a Graph}\label{graphmodule}

For a simple graph $G$ with $n$ vertices, 
let $M(G)$ be the $\ring$-submodule
of $\Mot$ 
generated by the $[A_G(s,r,k)]$.   
Let $M(G)_t$ be the submodule of $M(G)$ generated by the 
$[A_G(s,r,k)]$ for $s\leq t$.
Theorem~\ref{firstred} shows that 
$[A_G(s,r,k)]\in M(G)_k$.  Thus we have a finite filtration
$$
M(G)=M(G)_n\supset M(G)_{n-1} \supset \ldots \supset M(G)_{0}=\ring.
$$

The goals of this section are to compute the structure of $M(G)$ and
to show that, in fact, $M(G)\subset\Graphs_*$.  To do this we
introduce three special schemes: $K_G(s) = A_G(s,s,s)$,
$J_G(s)=\union_k A_G(s,s,k)$ and $H_G(s)= A_G(n,s,n)$.  Note that
$J_G(s)$ consists of the scheme of all pairs $(Q,f)\in A_G(s)$ with
$Q\in\Sym_s^s$; that is, there is no restriction on the rank of $f$.
Note also that $K_G(n)=H_G(n)$.

\begin{theorem}\label{th:Mgenerators}
\begin{enumerate} 
\item[(a)] $[A_G(s,r,k)]\in M(G)_d$ for $d=\min(s,r,k)$.
\item[(b)] $M(G)_t$ is spanned as an $\ring$-module by
the $[K_G(s)]$ for $s\leq t$.
\item[(c)] $M(G)_t$ is spanned as an $\ring$-module by
the $[J_G(s)]$ (resp. by $H_G(s)$) for $s\leq t$.
\end{enumerate}
\end{theorem}

\begin{proof} (a) and (b). Apply Theorem~\ref{firstred}
to obtain an expression for $[A_G(s,r,k)]$ as a $\Z[q]$-linear
combination of terms of the form $[A_G(k,j,k)]$ with
$j\leq\min(r,k)\leq s$.  Then apply Corollary~\ref{cor-secondred}
to obtain an expression for each $[A_G(k,j,k)]$ as a $\Z[q]$-linear
combination of terms of the form $[K_G(j)]$.

(c) To see that the $[J_G(s)]$ span, note that (a) implies that
$[J_G(s)]\equiv [K_G(s)]$ modulo $M(G)_{s-1}$.  To see that the the
$[H_G(s)]$ span, use the fact that $[H_G(s)]=\sigma [K_G(s)]$ for
$\sigma\in\calS$, a consequence of Corollary~\ref{cor-secondred}.
\end{proof}

Our interest in the $H_G(s)$ is based on the following lemma,  
which allows us to compare the $[H_G(s)]$ to the $[H_{DG}(s)]$.
The lemma is essentially a translation of Theorem 5.1 of \cite{stanley}
into our language.   As there are two graphs involved in the lemma, 
we write $n_G$ for the cardinality of $V(G)$.

\begin{lemma}\label{yuck} For $r\leq n_G+1$, 
\begin{equation}\label{eq:yuck}
[H_{DG}(r)] = a_G(r)[H_G(r)] + b_G(r)[H_G(r-1)] + 
c_G(r)[H_G(r-2)]
\end{equation}
with 
\begin{eqnarray*}
a_G(r) &=& q^{n_G+r}(q^{n_G+1} -1)\\
b_G(r) &=& q^{n_G+r-1}(q^{n_G+1} -1)(q-1)\\
c_G(r) &=& q^{n_G}(q^{n_G+1} -1)(q^{n_G+1} - q^{r-1})
\end{eqnarray*}
all polynomials in $\calS$.  
\end{lemma}
\begin{proof} 
By the Theorem~\ref{Dreduction}, 
$$
[A_{DG} (n_G+1,r,n_G+1)] =  (q^{n_G+1} - q^{n_G})[A_G(n_G+1, r, n_G)].
$$
Now applying Theorem~\ref{firstred} to 
$[A_G(n_G+1,r,n_G)]$ and expanding out
$[\Gr(n_G,n_G+1)]$ in terms of $q$ gives the result.

The polynomials $a_G, b_G$ and $c_G$ in the theorem are clearly in $\calS$
as long as they are nonzero. Inspection shows that this is the case under
the assumption that $r\leq n_G+1$. 
\end{proof}

There is a simpler identity relating $J_G(s)$ to $J_{DG}(s)$.

\begin{proposition}\label{pr:Jyuck} $[J_{DG}(s)]=q^s[J_G(s)]$.  
\end{proposition}
\begin{proof} The obvious map $J_{DG}(s)(\F_q)\to J_G(s)(\F_q)$ 
restricting $f$ from $V(DG)$ to $V(G)$ has fiber $\F_q^s$.
\end{proof}

A direct consequence of Proposition~\ref{pr:Jyuck} and 
Theorem~\ref{th:Mgenerators} (c) is the following.

\begin{theorem}\label{th:MDGMG} $M(DG)=M(G)$.
\end{theorem}

We are now ready to prove the main theorem of this section.  

\begin{theorem}\label{thm:graphmodulemain} 
$M(G)$ is equal to the $\ring$-module 
spanned by the functions $[Z_{D^kG}]$ for $k\geq 0$.  In particular,
$M(G)\subset\Graphs_*$.
\end{theorem}
\begin{proof} For the proof, let $N(G)$ be the $\ring$-module 
spanned by the functions $[Z_{D^kG}]$ for $k\geq 0$.  
Since $[K_G(n_G)]=[Z_G][\Gl_{n_G}]$, $[Z_G]\in M(G)$.  Thus
it follows from Theorem~\ref{th:MDGMG} that $N(G)\subset M(G)$.
To prove that $M(G)\subset N(G)$ we use Lemma~\ref{yuck} and an inductive
argument.  

By Theorem~\ref{th:Mgenerators} (c), it will be enough to show that 
$[H_G(s)]\in N(G)$ for all $s$.  Since $[H_G(n_G)]=[K_G(n_G)]$ this is 
obvious for $s=n_G$.  Now by Lemma~\ref{yuck} 
\begin{equation}
[H_{DG}(n_G + 1)] = b_G [H_G(n_G)] + c_G[H_G(n_G -1)]
\end{equation}
with $b_G, c_G\in\calS$.  (The first term on the right hand side of
(\ref{eq:yuck}) vanishes because $H_G(n_G+1)$ is empty).  We know that
$[H_{DG}(n_G + 1)]$ and $[H_G(n_G)]$ are in $N(G)$.  Thus $[H_G(n_G
-1)]\in\calS$.

We then assume inductively that $[H_G(n_G -i)]\in N(G)$ for all $i\leq
a$ and for all graphs $G$.  Another application of Lemma~\ref{yuck}
shows us that
\begin{equation}
\begin{split}
[H_{DG}(n_G -(a-1))] = &a_G[H_G(n_G - (a-1))] + b_G[H_G(n_G -a)]\\ 
                       &+c_G[H_G(n_G - (a+1))].
\end{split}
\end{equation}
By induction, the left-hand side and the two first terms on the right
hand side are in $N(G)$.  Thus, as $c_G\in\calS$, 
$[H_G(n_G - (a+1))]\in N(G)$ as well.
\end{proof}

\section[matroids]{Matroid Theory}\label{sec:matroidtheory}

A matroid $M$ consists of a finite set $E$ called the edges of the matroid
and a rank function $\rho:2^E\to\N$ satisfying the following axioms
\begin{enumerate}
\item $\rho(E)\leq \num{E}$.
\item For $X\subset Y\subset E$, $\rho (X)\leq \rho(Y)$.
\item For any $X,Y\subset E$, 
\begin{equation}
\rho(X\cup Y) + \rho(X\cap Y) \leq \rho(X) + \rho(Y).
\end{equation}
\end{enumerate}
The integer $\rho(E)$ is said to be the {\em rank} of the matroid.

Matroids were introduced by H. Whitney \cite{whitney} as a simultaneous
generalization of matrices and graphs.  An excellent modern reference
for matroid theory is \cite{oxley}.  

\subsubsection{Representability}
A matroid $M$ of rank $r$ is said to be representable over a 
field $K$ if there is a 
function $f:E\to K^r$ such that the dimension of the span of the set 
$f(X)$ is equal to $\rho(X)$ for all $X\subset E$.

\subsubsection{Matroids from matrices}\label{subsubsec:mfromm}
To every subset  $E\subset K^s$ there is naturally a matroid
$M$ representable over $K$ given by setting $\rho(X)=\dim\spn{X}$
for every $X\subset E$.  

Let $K^{n+1}-\{0\}\stackrel{\pi}{\mapsto}\bP^n(K)$ be the natural map
taking a nonzero vector in $v\in K^{n+1}$ to the line $Kv$.  Suppose
$E\subset\bP^n(K)$.  Then any set-theoretic splitting
$\sigma:\bP^n(K)\mapsto K^{n+1}$ gives a subset $\sigma(E)$ of
$K^{n+1}$ and, thus, defines a matroid.  It is easy to see that this
matroid is independent of the splitting $\sigma$.  Thus, since such
splitting always exist, $E$ defines a matroid.

\subsubsection{Representation schemes}
For any matroid $M$ and a locally free sheaf $W$ over a base $S$, let
$X(M,W)$ be the subscheme of $\Fun(E,W)$ consisting of all $f$ whose
restrictions to $\Fun(X,W)$ lie in $\Fun_{\rho(X)}(X,W)$ for all
$X\subset E$.  This is the scheme of representations of $M$ in $W$.
When $S=\Spec\Z$ and $W=\calO_S^s$, we write $X(M,s)$ for $X(M,W)$ as
in the introduction.  For a field $K$, $X(M,s)(K)$ is the set of all
maps $f:E\to K^s$ such that $\dim\spn{f(X)}=\rho(X)$ for all $X\subset
E$.  That is, $X(M,s)(K)$ is the set of all representations of $M$ in
$K^s$.  When $r$ is the rank of $M$, we write $X(M)$ for $X(M,r)$.
$X(M)(K)$ is non-empty if and only if $M$ is representable over $K$.
Clearly $X(M,s)(K)$ is non-empty if and only if $s>r$ and $X(M)(K)$ is
non-empty.

\begin{definition}
Let $\Matroids$ be the $\ring$-module generated by all functions
of the form $[X(M)]$.
\end{definition}
\begin{remark} It is easy to see that $[X(M,s)]=[\Gr(r,s)][X(M)]$.
Thus $\Matroids$ is the same as the $\ring$-module generated by
all functions of the form $[X(M,s)]$.
\end{remark}

In the next section we show that $\Matroids\subset\Graphs_*$.   

\section[counterexample]{A Counterexample to Kontsevich's Conjecture}
\label{sec:counterexample}

Let $G$ be a graph, $V$ the set of its vertices, $U\subset 2^V$.  A
function $\pi:U\mapsto \N$ will will be called a {\em partially
defined rank function for $V$}.  Notice that the data of a partially
defined rank function $\pi$ determines $U=\dom(\pi)$.  Associated to
every such function we have a scheme defined as follows:

\begin{definition}\label{maindefp}  
$J_G(s,\pi)$ is the scheme of all 
of all $(Q,f)\in J_G(s)$ such that $f_{|H}$ has rank $\rho(H)$
for all $H\in\dom{\pi}$. 
\end{definition}

\begin{theorem}\label{pi} 
For every $G$ and every partially defined rank function
$\pi$ for $V(G)$, $[J(s,\pi)]\in\Graphs_*$.
\end{theorem}
\begin{proof}

The proof is by induction on the cardinality of $\dom(\pi)$.
If $\dom(\pi)$ is empty, $J_G(s,\pi)=J_G(s)$.  Thus the result
follows from Theorem~\ref{thm:graphmodulemain}.

Now assume the result holds for all graphs $G$ and all $\pi$ such that
$\num{\dom{\pi}}\leq a$. Let $W\subset 2^V$ be a set of subsets with
$a+1$ elements, let $H\in W$ and let $U=W-\{H\}$.  Let $\pi:U\to\N$ be
a partially defined rank function, and let $\pi_i:W\to\N$ be the
extension of $\pi$ to $W$ such that $\pi_i(H)=s-i$.  Clearly, any
partially defined rank function with domain $W$ is of the form $\pi_i$
for some $\pi:U\to\N$ and some $i\in[0,s]$.

Now for each $t\in\N$ we define a graph $G_t$ as follows:
$G_t$ is the graph obtained from $G$ by adjoining $t$ disjoint
vertices $y_1,\ldots, y_t$ and connecting each of the $y_i$ by 
edges only to the vertices in $H$.  Thus $V(G_t) = V(G)\union Y$
where $Y=\{y_1,\ldots, y_t\}$, and 
$$
E(G_t)=E(G)\union \{e_{hy}\}_{h\in H\atop y\in Y}.
$$

Since $V(G)\subset V(G_t)$, $U\subset 2^{V(G_t)}$.   Thus we can
consider $\pi$ as a partially defined rank function for $V(G_t)$.

The result will follow from the following equation:
\begin{equation}
[J_{G_t}(s, \pi)] = \sum_{i=0}^s q^{ti} [J_G(s, \pi_i)]
\end{equation}

To see that the equation holds,  note that we can stratify the $\F_q$ 
points of $J_{G_t}(s,\pi)$
according to the dimension of the span of $f(H)$.  
Let $J_{G_t}(s,\pi)_i$ be the stratum where this dimension is $s-i$.
This stratum maps to $J_G(s,\pi_i)$ by restricting $f$ from 
$V(G_t)$ to $V(G)$.  The fiber of map
above any point $(Q,f)$ is an affine space $\A^{ti}$.  
This is because the only condition on the $f(y_i)$ is that they 
be orthogonal to the span of $f(H)$.  Thus, as the bilinear
form $Q$ is always nondegenerate, they must lie in a linear
subspace of dimension $i$.  

To complete the proof, note that by varying the $t$ from
$0$ to $s$ we obtain
a system of equations for the $[J_G(s, \pi_i)]$ in terms of the 
$[J_{G_t}(s, \pi)]$.  Solving this system for the $J_G(s,\pi_i)$
using Cramer's rule, we have to invert a Vandermonde determinant
which lies in $\calS$.
Thus, as we assumed by induction that $[J_{G_t}(s,\pi)]$ lies
in $\Graphs_*$, it follows that each $[J_G(s,\pi_i)]$ lies
in $\Graphs_*$ as well.
\end{proof}

Let G be a discrete graph (that is $E(G)$ is empty). 
In this case,
if  $\pi$ is a partially defined rank function then
\begin{equation}
[J(s,\pi)]=[\Sym^s_s][L(s,\pi)]
\end{equation}
where $L(s,\pi)$ is the scheme consisting of all
$f\in\Fun(V,\calO_{\Spec\Z}^s)$ such that $f$ restricts to
$\Fun_{\pi(H)}(H,\calO_{\Spec\Z}^s)$ for all $H\in\dom(\pi)$.  To see
this, note that the definition of $J(s,\pi)$ makes it clear that Q
does not enter in the definition of the J's for the discrete graph.
And only the vertex set $V$ is needed for the definition of the $L$'s
since $G$ is discrete.

As $\Sym^s_s\in\calS$, it follows that the $L$'s are all in
$\Graphs_*$.  Now note that, if $M$ is a matroid, with rank function
$\rho:2^E\to\N$, then $X(M,s)=L(s,\rho)$.  Thus we obtain the
following theorem:
\begin{theorem}\label{thm:ming} $\Matroids\subset\Graphs_*$.
\end{theorem}

It is now possible to see directly that
Conjecture~\ref{conjSt} and thus 
Conjecture~\ref{conjK}
are false.   Let $M$ be the Fano matroid.  This is 
a rank $3$ matroid whose edge set $E$ is the set $\bP^2(\F_2)$.    
This matroid is representable
over a field $\F_q$ if and only if $2|q$ (see \cite{welsh} 
Chapter 9).   

Thus the function $[X(M)]$ is supported on the set of $q$ such
that $2|q$.  It follows that $[X(M)]$ can not be a rational
function.  And this contradicts Conjecture~\ref{conjSt}
by Theorem~\ref{thm:ming}.

\section[Mnev]{Representation problem of Matroids}\label{sec:mnev}

Our objective in this section is to show that $\Matroids=\Mot$
and thus that 
$\Graphs_* = \Mot.$
This will follow from the known results on the Matroid
representation problem.

We saw in the previous section that $[L(k,\pi)]$ were in $\Graphs_{*}$
even if $\pi$ is only partially defined. It suffices
therefore to show that the $\ring$-module generated by all
functions of the form $[L(k,\pi)]$ is all of $\Mot^+$.  This was in
essence proved by Mn\"ev \cite{Mnev, mnevb} as the unoriented matroid
component of a more difficult theorem concerning the representation
spaces of oriented matroids (see also \cite{gunzel, shor}).  It was
independently proved by Bokowski and Sturmfels \cite{BokowskiSturmfels,
sturmfelsa}.  Moreover, the idea of the proof using von Staudt's
``algebra of throws'' goes back at least to \cite{maclane} (see
\cite{kung} for an enlightening explication).  However, as we have
been unable to extract a proof of the exact statement we need from the
literature, we give a sketch of the proof in our context.

\begin{theorem}[Mn\"ev,Sturmfels]\label{th-MS}
Let $X$ be a quasi-projective scheme of finite type over $\Z$, then
there is a set $V$, a set of subsets $W$ of $V$, a function $\pi :W
\mapsto \Z$, and an element $\sigma\in\calS$ so that
$$\sigma [X] = [L(3,\pi)].$$
\end{theorem}

\begin{remark}
 
\begin{enumerate}

\item The theorems in Matroid theory are not in such a direct form
because, in Matroid theory we are committed to declare the rank of all
the subsets of $V$. Our partially defined $\pi$ does not have this
problem. By inclusion-exclusion principles the $\ring$-module
generated by all functions of the form $[L(k,\pi)]$ where $\pi$ may
only be partially defined is same as the $\ring$-module
generated by all functions of the form $[L(k,\pi)]$ where $\pi$ is
defined on all subsets of $V$.

\item Note that any scheme of finite type/$\Z$ is a finite disjoint union
of quasiprojective schemes/$\Z$.
\end{enumerate}

\end{remark}
\begin{proof} 
The proof follows essentially from the following observations

\begin{enumerate}
\item $4$ elements in $P^2$ such that any $3$ are linearly independent
can by a unique automorphism of $P^2$ in PGL(2), be assumed to be
$(1,0,0), (0,1,0), (0,0,1)$ and $(1,1,1)$.
\item if given two points $(x,0,1)$ and $(x',0,1)$ on the $X-axis$,
then by drawing lines alone through the 4 points above and these two
points, we can locate $(x+x',0,1)$,$(xx',0,1)$,$(-x,0,1)$. Finding
intersection of lines can be translated as a vector which lies on both
lines, and hence as a condition on linear dependence. These
constructions can be found for example in the proof of Theorem 2.2 of
\cite{BokowskiSturmfels} .
\item Iterating these constructions, given $(x_1,x_2,\dots,x_n)$ we
can determine the points $(f(x_1,x_2,\dots,x_n),0,1)$, where f is a
polynomial with integer coefficients by just drawing lines starting
from the configuration of the four given points and the points
$(x_i,0,1)$. Setting $f(x_1,\dots,x_n)$ either equal to zero or not
equal to zero is just another spanning condition: A condition on
whether
$$(f(x_1,x_2,\dots,x_n),0,1),(0,0,1)$$ 
is linearly independent
or not.
\item The cone over any quasi-projective scheme/$\Z$ can be written as
a set of equalities and a set of nonequalities in a finite set of
variables $(x_1,\dots,x_n)$. Note that we can also have conditions of
the form $n=0$ in the list.
\end{enumerate}
\end{proof}

\section[forest]{Forests}\label{forests}

In this section we prove that $[J_G(s)]\in\Z[q]$ whenever $G$ is a
forest.  It follows that $M(G)=\Z[q]$ for such graphs.  To do so we
need to introduce a two operations on graphs.

Let $v\in V(G)$.  We obtain a graph $I_v(G)$ by adding one edge $e$
connected to $v$ and one new vertex $w$ connected to $e$.  That is, we
insert an edge at $v$.  Clearly, a graph is tree if and only if it can
be obtained from the graph with one vertex by successive applications
of $I_v$ for various $v$.  A graph is a forest if and only if it can
be obtained from the empty graph $\emptyset$ by successive
applications of $I_v$ and the operation $D$.  We write $D_n$ for the
graph $D^n\emptyset$.

We define $R_v$ to be the graph obtained from $G$ by deleting
$v$ and all edges meeting it.  Note that if $G$ is a forest
and $v$ is any vertex in $G$, $R_v G$ is also a forest.

\begin{theorem}\label{th-forests} 
Let $G$ be a graph with $v\in V(G)$.
\begin{enumerate}
\item[(a)] $[J_{D_n}(s)] = q^{ns}[\Sym^s_s]$.
\item[(b)] $[J_{DG}(s)]  = q^s J_G (s).$
\item[(c)] $[J_{I_v}(s)] = q^{s-1}\left(J_{G}(s) + (q-1) J_{R_vG}(s)\right)$
\end{enumerate}
\end{theorem}
\begin{proof} For (b) let $w$ be the vertex in $V(DG)-V(G)$.  
Then the map $J_{DG}(s) \to J_{G}(s)\times \F_q^s$ given by
$(Q,f)\mapsto (Q, f|_{V(G)}, f(w))$ is an isomorphism.
For (a) assume first that $n=0$.  Then tracing through the definitions
one sees that $J_{\emptyset}(s) = \Sym^s_s$.  The rest of (a) follows
by induction from (b).

For (c) we work over $\F_q$ and 
consider the map $\pi: J_{I_v(G)}(s) \to J_{G}(s)$ given
by $(Q,f)\mapsto (Q, f|_{V(G)})$.  The fiber of $\pi$ above a point
$(Q, g)\in J_{G}(s)$ depends on whether $g(v)$ is $0$ or not.
Let $J^0_G(s)$ (resp. $J^{\times}_G(s)$) be the set where $g(v)=0$
(resp. $g(v)\neq 0$).   Above a point $(Q,g)\in J^0_G(s)$ the fiber 
of $\pi$ will have $q^s$ points.   Above a point 
$(Q,g)\in J^{\times}_G(s)$ the fiber will have $q^{s-1}$ points
since $Q$ is non-degenerate.  Thus 
\begin{equation}
|J_{I_v(G)}(s)| = q^{s-1} |J^{\times}_G(s)| + q^s|J^0_G(s)|
\end{equation}
The result now follows from the observation that 
$|J^0_G(s)| = |J_{R_vG}(s))|$.
\end{proof}

\begin{corollary} For $F$ a forest, $[Z_F]\in\Z[q]$.
\end{corollary} 
\begin{proof}  An easy induction using Theorem~\ref{th-forests} 
shows that $[J_F(s)]\in\Z[q]$ for any $s$.  Thus $M(F)=\ring$.  It
follows from Theorem~\ref{thm:graphmodulemain} that $[Z_F]\in\ring$.
But this implies that $[Z_F]\in\Z[q]$ by
Proposition~\ref{th-prational}.
\end{proof}

The next corollary follows 
from Theorem~\ref{th-forests} and the results in 
Section~\ref{matrix-tree}.

\begin{corollary}\label{co-forests}
Let $F$ be a forest with $r$ vertices contained
in a complete graph $K_s$.  Let $G= K_s -F$. 
\begin{enumerate} 
\item $[Z^o_G]\in\Z[q]$.
\item If $s>r$, then $[X_G]\in\Z[q]$.
\end{enumerate}
\end{corollary}

\bibliographystyle{plain}
\def\noopsort#1{}

\end{document}